\documentclass[12pt,a4paper,reqno]{amsart}
\usepackage{graphicx}
\usepackage[all]{xy}
\usepackage{amssymb}
\usepackage{amsmath}
\usepackage[utf8]{inputenc}
\usepackage{amsthm}

\newtheorem{proposition}{Proposition}
\newtheorem{definition}{Definition}
\newtheorem{theorem}{Theorem}
\newtheorem{example}{Example}

\newcommand{\C}{\mathbf{C}}

\newcommand{\Aut}{\mathrm{Aut}}

\begin{document}

\title[Internal groupoids as involutive-2-links]{Internal groupoids as involutive-2-links}


\author{Nelson Martins-Ferreira}
\address[Nelson Martins-Ferreira]{Instituto Politécnico de Leiria, Leiria, Portugal}
\thanks{ }
\email{martins.ferreira@ipleiria.pt}

\begin{abstract}
Regardless of its environment, the category of internal groupoids is shown to be equivalent to the full subcategory of involutive-2-links that are unital and associative. The new notion of involutive-2-link originates from the study of triangulated surfaces and their application in additive manufacturing and 3d-printing. Thus, this result establishes a bridge between the structure of an internal groupoid and an abstract triangulated surface. An example is provided which can be thought of as a crossed-module of magmas rather than groups.



\end{abstract}

\maketitle

\date{Received: date / Accepted: date}

The purpose of this note is to build a bridge betwwen the study of internal groupoids  and the study of triangulated surfaces. The structure of an abstract triangulated surface, as described in \cite{NMF23}, has motivated the search for an analogous model to an internal groupoid. The result is presented here under the name \emph{involutive-2-link} with its two main properties: unitary and associativity.

\begin{theorem}\label{thm:one}
Let $\C$ be any category. The category of internal groupoids is equivalent to the full subcategory of unital and associative involutive-2-links.
\end{theorem}

An \emph{involutive-2-link} is a morphism $m\colon{A\to B}$ equipped with two interlinked involutions on its domain. More precisely, it consists of a triple $(\theta,\varphi,m\colon{A\to B})$ with $\theta,\varphi\colon{A\to A}$ such that $\theta^2=\varphi^2=1_A$ and $\theta\varphi\theta=\varphi\theta\varphi$. 
Note that the subgroup of $\Aut(A)$, generated by $\theta$ and $\varphi$, is the dihedral group of order 6.

A morphism between involutive-2-links, say from $(\theta,\varphi,m\colon{A\to B})$ to $(\theta',\varphi',m'\colon{A'\to B'})$ is a pair of morphisms $(f\colon{A\to A'},g\colon{B\to B'})$ such that $f\theta=\theta' f$, $f\varphi=\varphi' f$ and $m'f=gm$.

\begin{definition}\label{def:1}
Let $\C$ be any category. An \emph{involutive-2-link} structure in $\C$, say $(\theta,\varphi,m\colon C_2\to C_1)$, is said to be:
\begin{enumerate}
\item \emph{unital} when the two pairs of morphisms $(m,m\theta)$, $(m,m\varphi)$ are jointly monomorphic and there exist morphisms $e_1,e_2\colon C_1\to C_2$ such that
\begin{align}
    me_1=1_{C_1}=me_2\label{eq: main 1 cp}\\
    \theta e_2=e_2,\quad \varphi e_1=e_1\label{eq: main 2 cp}\\
    m\theta\varphi e_2=m\varphi\theta e_1\label{eq: main 3 =i cp}\\
    m\theta e_1 m\varphi=m\varphi e_2 m\theta\label{eq: main 4 ed=ec cp}\\
    m\theta e_1 m=m\theta e_1 m\theta \label{eq: main 5 contract e1} \\
    \quad m\varphi e_2 m=m\varphi e_2 m\varphi.\label{eq: main 5 contract e2} 
\end{align}

\item \emph{associative} when the pair $(m\varphi,m\theta)$ is bi-exact (see diagram (\ref{diag: di-graph pi_1 pi_2}) bellow with $m\varphi$ as $\pi_1$ and $m\theta$ as $\pi_2$) and the induced morphisms $m_1,m_2\colon C_3\to C_2$, determined by (see diagram (\ref{diag: zig-zag completed}))
\begin{align*}
    \pi_1m_1=mp_1,\quad \pi_2 m_1=\pi_2 p_2\\
    \pi_1 m_2=\pi_1 p_1, \quad \pi_2 m_2=m p_2
\end{align*}
are such that $mm_1=mm_2$.
\end{enumerate}
\end{definition}

A pair of parallel morphisms (or a digraph) is said to be  \emph{bi-exact} if when considered as a span it can be completed into a commutative square which is both a pullback and a pushout and moreover, if considered as a cospan, it can be completed into another commutative square which is both a pullback and a pushout. In other words, a digraph such as 
\begin{equation}\label{diag: di-graph pi_1 pi_2}
\xymatrix{C_2 \ar@<0.5ex>[r]^{\pi_2}\ar@<-0.5ex>[r]_{\pi_1} & C_1}
\end{equation}
is \emph{bi-exact} precisely when the zig-zag 
\begin{equation}\label{diag: zig-zag}
\xymatrix{
& C_2\ar[d]^{\pi_1}\\ C_2 \ar[r]_{\pi_2}\ar[d]_{\pi_1} & C_1\\ C_1 }
\end{equation}
can be completed with two commutative squares 
\begin{equation}\label{diag: zig-zag completed}
\xymatrix{ C_3\ar@{-->}[r]^{p_2} \ar@{-->}[d]_{p_1}
& C_2\ar[d]^{\pi_1}\\ C_2 \ar[r]_{\pi_2}\ar[d]_{\pi_1} & C_1 \ar@{-->}[d]^{c}\\ C_1 \ar@{-->}[r]^{d} & C_0 }
\end{equation}
which  are both simultaneously a pullback and pushout. Such squares are also called exact squares, bicartesian squares, Dolittle diagrams or pulation squares \cite{Ban}.

The notion of a bi-exact pair of parallel morphisms is a way to study internal groupoids in arbitrary categories, even thought pullbacks may not be available as canonical constructions. However, since the results are invariant via the Yoneda embedding, our proofs are based in the ambient category of sets and maps. Nevertheless,  details are given as if working in a context where pullbacks have to be considered as a property of commutative squares.

The functor $F$ from the category of internal groupoids to the category of involutive-2-links is defined via the assignment
\begin{equation}\label{diag: internal groupoid assignment}\xymatrix{C_2 \ar@<1.5ex>[r]^{\pi_2}\ar@<-1.5ex>[r]_{\pi_1}\ar[r]|{m} & C_1 \ar@(ur,ul)[]_{i}  \ar@<1ex>[r]^{d} \ar@<-1ex>[r]_{c} & C_0 \ar[l]|{e} }
\mapsto
\xymatrix{C_2 \ar@(ul,dl)[]_{\theta,\varphi} \ar[r]^{m} & C_1},
\end{equation}
with $\theta=\langle i\pi_1,m\rangle$, $\varphi=\langle m,i\pi_2\rangle$ and it is full and faithful.  
Indeed, let us consider an internal groupoid (\cite{NMF22a}, see also \cite{BorceuxJanelidze}, Section 7.1) as a diagram of the form
\begin{equation}\label{diag: internal groupoid}\xymatrix{C_2 \ar@<1.5ex>[r]^{\pi_2}\ar@<-1.5ex>[r]_{\pi_1}\ar[r]|{m} & C_1 \ar@(ur,ul)[]_{i}  \ar@<1ex>[r]^{d} \ar@<-1ex>[r]_{c} & C_0 \ar[l]|{e} }
\end{equation}
such that
\begin{eqnarray}
de = 1_{C_1} =ce \\
d m = d \pi_2,\quad c m = c \pi_1,\quad d\pi_1=c\pi_2\\
di=c,\quad ci=d, \quad i^2=1_{C_1},\quad ie=e
\end{eqnarray}
and satisfying the following further properties:
\begin{enumerate}

\item[(a)] the commutative square 
\begin{equation}\label{diag: d c square}
\vcenter{
\xymatrix{C_2 \ar[r]^{\pi_2}\ar[d]_{\pi_1} & C_1 \ar[d]^{c}\\
 C_1\ar[r]^{d} & C_0}}
\end{equation}
is a pullback square;

\item[(b)] $m\langle 1_{C_1},ed\rangle=1_{C_1}=m\langle ec,1_{C_1}\rangle$;

\item[(c)] $m\langle 1_{C_1},i\rangle=ec,\quad m\langle i,1_{C_1}\rangle=ed$;\label{property c}

\item[(d)] the cospan $\xymatrix{C_2\ar[r]^{d\pi_2} & C_0 & C_1 \ar[l]_{c}}$ can be completed into a pullback square
\begin{equation}\label{diag: dpi_2 c square}
\vcenter{
\xymatrix{C_3 \ar[r]^{p_2}\ar[d]_{p_1} & C_1 \ar[d]^{c}\\
 C_2\ar[r]^{d\pi_2} & C_0}}
\end{equation}

\item[(e)] $m(1\times m)=m(m\times 1)$, where $(1\times m),(m\times 1)\colon C_3\to C_2$ are morphisms uniquely determined as
\begin{align*}
    \pi_2(m\times 1)=p_2\\
    \pi_1(m\times 1)=m\\
    \pi_2(1\times m)=m\langle \pi_2p_1,p_2 \rangle\\
    \pi_1(1\times m)=\pi_2.
\end{align*}
\end{enumerate}

The functor $F$ takes an internal groupoid such as (\ref{diag: internal groupoid}), forgets the underlying reflexive graph 
\begin{equation}\label{diag: underlying reflexive graph}\xymatrix{C_1   \ar@<1ex>[r]^{d} \ar@<-1ex>[r]_{c} & C_0 \ar[l]|{e} },
\end{equation}
keeps the morphism $m\colon{C_2\to C_1}$ (the multiplicative structure of the internal groupoid) and contracts the remaining information as two endomorphisms $\theta,\varphi\colon{C_2\to C_2}$ of the form
\begin{equation}\label{eg: theta and varphi def}
 \theta=\langle i\pi_1,m\rangle \quad \varphi=\langle m,i\pi_2\rangle.
\end{equation}
As a consequence, we have
\begin{align}
    m\varphi=\pi_1,\quad m\theta=\pi_2\label{eq: theta and varphi 1}\\
    \pi_1\varphi= m,\quad \pi_1\theta=i\pi_1\label{eq: theta and varphi 2}\\
    \pi_2\varphi=i\pi_2,\quad \pi_2\theta= m.\label{eq: theta and varphi 3}
\end{align}
 The conditions $\theta^2=\varphi^2=1_{C_2}$ and $\theta\varphi\theta=\varphi\theta\varphi$ are easily verified. Hence, the functor is well defined and it is clearly faithful. 

In order to see that the functor $F$ is full, let us consider two internal groupoids, say
\begin{equation}\label{diag: internal groupoid C}\xymatrix{C_2 \ar@<1.5ex>[r]^{\pi_2}\ar@<-1.5ex>[r]_{\pi_1}\ar[r]|{m} & C_1 \ar@(ur,ul)[]_{i}  \ar@<1ex>[r]^{d} \ar@<-1ex>[r]_{c} & C_0 \ar[l]|{e} }
\end{equation}
and 
\begin{equation}\label{diag: internal groupoid C'}\xymatrix{C'_2 \ar@<1.5ex>[r]^{\pi'_2}\ar@<-1.5ex>[r]_{\pi'_1}\ar[r]|{m'} & C'_1 \ar@(dr,dl)[]^{i'}  \ar@<1ex>[r]^{d'} \ar@<-1ex>[r]_{c'} & C'_0 \ar[l]|{e'} }
\end{equation}
denoted respectively by $C$ and $C'$. Let us assume the existence of a morphism of involutive-2-links from $F(C)$ to $F(C')$, that is, a pair of morphisms $f_i\colon{C_i\to C'_i}$, with $i=1,2$ such that $\theta'f_2=f_2\theta$, $\varphi'f_2=f_2\varphi$ and $m'f_2=f_1m$, with $\theta,\varphi,\theta',\varphi'$ the respective involutions associated with $F(C)$ and $F(C')$. We need to show that the pair $(f_2,f_1)$ can be extended to a morphism of internal groupoids 
\begin{equation}\label{diag: internal groupoid C to C'}\xymatrix{C_2\ar[d]_{f_2} \ar@<1.5ex>[r]^{\pi_2}\ar@<-1.5ex>[r]_{\pi_1}\ar[r]|{m} & C_1 \ar[d]_{f_1}\ar@(ur,ul)[]_{i}  \ar@<1ex>[r]^{d} \ar@<-1ex>[r]_{c} & C_0 \ar@{-->}[d]^{f_0}\ar[l]|{e} \\
C'_2 \ar@<1.5ex>[r]^{\pi'_2}\ar@<-1.5ex>[r]_{\pi'_1}\ar[r]|{m'} & C'_1 \ar@(dr,dl)[]^{i'}  \ar@<1ex>[r]^{d'} \ar@<-1ex>[r]_{c'} & C'_0 \ar[l]|{e'}. }
\end{equation}
First observe that $f_2(x,y)=(f_1(x),f_1(y))$ since $\pi'_1f_2=m'\varphi'f_2=m'f_2\varphi=f_1m\varphi=f_1\pi_1$ and similarly $\pi'_2f_2=f_1\pi_2$. This means that the hypotheses $\theta'f_2=f_2\theta$, $\varphi'f_2=f_2\varphi$ and $m'f_2=f_1m$ are translated, respectively, as
\begin{align*}
(f_1(x)^{-1},f_1(x)f_1(y))&=(f_1(x^{-1}),f_1(xy))\\
(f_1(x)f_1(y),f_1(y)^{-1})&=(f_1(xy),f_1(y^{-1}))\\
f_1(x)f_1(y)&=f_1(xy)
\end{align*}
from which we conclude $i'f_1=f_1i$. We also have $f_1ed(x))=f_1(x^{-1}x)=f_1(x^{-1})f_1(x)=f_1(x)^{-1}f_1(x)=e'd'f_1(x)$ and  $f_1ec(x)=e'c'f_1(x)$, which give
\begin{align*}
\langle 1,e'd'\rangle f_1=f_2\langle 1,ed\rangle\\
\langle 1,e'c'\rangle f_1=f_2\langle 1,ec\rangle
\end{align*}  
and permits the definition of $f_0$ either as $d'f_1e$ or as $c'f_1e$. Hence, the triple $(f_2,f_1,f_0)$ is a morphism of internal groupoids from $C$ to $C'$, showing that the functor $F$ is full.

Let us observe that even when $i$ is not made explicit, the two involutions $\theta$ and $\varphi$ are still uniquely determined because the object  $C_2$ can be presented not only as the pullback of $d$ and $c$ but also as the kernel pair of $d$ or the kernel pair of $c$ and hence both pairs $(m,\pi_1)$ and $(m,\pi_2)$ are in particular jointly monomorphic.    This fact suggests the possibility of considering an even simpler structure to describe internal groupoids by using only $\theta$ or $\varphi$ together with the multiplication $m$. However, this would give rise to a different structure which requires further investigation. Nevertheless, it is possible that the bridge with triangulated structures  \cite{NMF23} will be widened by the new structure to be found.
Moreover, the common denominator to the unitary and associativity properties is the requirement that the two pairs $(m,m\theta)$ and $(m,m\varphi)$ are  jointly monomorphic which reinforces the possibility of having, say, the existence of $\varphi$ as a property of the pair $(m,m\theta)$. 


\vspace*{1cm}

In order to prove Theorem \ref{thm:one} it is readily seen that that if $(\theta,\varphi,m)$ is obtained from an internal groupoid by applying the functor $F$ then it is a unital and associative involutive-2-link. On the other hand, if $(\theta,\varphi,m)$ is a unital and associative involutive-2-link, then the fact that the pairs $(m,m\theta)$ and $(m,m\varphi)$ are jointly monomorphic uniquely determines the morphisms $e_1$ and $e_2$ which are required to exist by the unitary property and hence fulfill the properies (b) and (c) of an internal groupoid. Indeed, the morphism $i\colon C_1\to C_1$ is obtained by condition (\ref{eq: main 3 =i cp}) either as $i=m\theta\varphi e_2$ or as $i=m\varphi\theta e_1$. The morphism $e\colon C_0\to C_1$ is uniquely determined by condition (\ref{eq: main 4 ed=ec cp}) as such that $ed=m\theta e_1$ and $ec=m\varphi e_2$ where $d$ and $c$ are obtained as in diagram $(\ref{diag: zig-zag completed})$ (with $m\varphi$ as $\pi_1$ and  $m\theta$ as $\pi_2$, which consequenctly also gives the property (a) of an internal groupoid because the pair $(m\theta,m\varphi)$ is bi-exact).  Conditions (\ref{eq: main 1 cp}), (\ref{eq: main 5 contract e1}) and (\ref{eq: main 5 contract e2})   assert the contractibility of the pairs $(m,m\theta)$ and $(m,m\varphi)$ in the sense of Beck (see \cite{MacLane}, p. 150). Condition (\ref{eq: main 2 cp}) is a central ingredient and gives $e_1e=e_2e$ from which the conditions $dm=d\pi_2$ and $cm=c\pi_1$   are deduced, thus permitting to define the two morphisms $m_1$ and $m_2$ from the fact that the pair $(m\theta,m\varphi)$ is bi-exact.

\vspace*{.5cm}
 The remaining details in the proof are easily obtained. Let us turn our attention to an example that can be seen as a generalization of crossed-modules from groups to magmas \cite{NMF22b,NMF23b}.

\begin{example}\label{e.g.:1}
Let $X=(X,\cdot)$ be a magma with a distinguished element $1\in X$ and $B$ be a non-empty set, with $0\in B$, together with maps $\bar{()}\colon{X\to X}$, $f\colon{X\times B\times X\times B\to X}$ and $g\colon{X\times B\to B}$ such that 
\begin{eqnarray}
f(y\cdot x,b,y'\cdot x',b')=f(y,g(x,b),y',g(x',b'))\cdot f(x,b,x',b')\label{eq:1 ff}\\
 g(y\cdot x,b)=g(y,g(x,b))\label{eq:2 gg}\\
 g(1,0)=0,\quad 1\cdot 1=1\label{eq:3 g(0,1) and 11=1}\\
 g(\bar{x}(\bar{y}(yx)),b)=b,\quad g(\bar{x}x,b)=b
 .\label{eq:4 abc}
\end{eqnarray}
Note that sometimes $x\cdot y$ is written as $x y$. Let us also consider the sets
\begin{eqnarray}
C_1=\{(x,b)\in X\times B\mid f(x,b,0,0)=x=f(0,0,x,b)\}\\
C_2=\{(y,x,b)\in X^2\times B\mid (y,g(x,b)),(x,b)\in C_1\},
\end{eqnarray}
and the formulas
\begin{eqnarray}
m(y,x,b)=(yx,b)\\
\theta(y,x,b)=(\bar{y},yx,b)\\
\varphi(y,x,b)=(yx,\bar{x},g(\bar{y}(y x),b)).
\end{eqnarray}
\end{example}
The following propositions refer to the structure of Example~\ref{e.g.:1} and should be considered as simple observations.

\begin{proposition}\label{prop:1} The maps $\theta,\varphi\colon{X^2\times B\to X^2\times B}$ are involutions if and only if the conditions
\begin{equation}\label{eq:5a tfae}
\bar{\bar{x}}=x,\quad \bar{y}(yx)=x,\quad (yx)\bar{x}=y
\end{equation}
hold for all $x,y\in X$. In addition, the further condition $\theta\varphi\theta=\varphi\theta\varphi$ is satisfied if and only if the two extra conditions 
\begin{equation}\label{eq:5b tfae}
 x(\overline{yx})=\bar{y},\quad (\overline{yx})y=\bar{x}
\end{equation}
are also satisfied for all $x,y\in X$. 
\end{proposition}

Let us restrict our attention to the subsets $C_1$ and $C_2$.

\begin{proposition}\label{prop:1.5}
If $(y,x,b)\in C_2$ then $(yx,b)\in C_1$.
\end{proposition}

Hence, $m(y,x,b)=(yx,b)$ is a well defined map $m\colon C_2\to C_1$.

\begin{proposition}\label{prop:2}
The formulas $\theta$, $\varphi$ are well defined maps $C_2\to C_2$ if and only if the following condition holds:
\begin{equation}\label{eq:6 tfae}
\text{if $(y,x,b)\in C_2$ then $(\bar{x},\bar{y},g(yx,b)\in C_2$.}
\end{equation}
\end{proposition}

Let us give sufficient conditions for the structure $(\theta,\varphi,m\colon C_2 \to C_1)$ to be (a well defined)  involutive-2-link. 

\begin{proposition}\label{prop:3} The structure $(\theta,\varphi,m\colon C_2 \to C_1)$ is (a well defined) involutive-2-link as soon as the following two conditions hold:
\begin{eqnarray}
\text{if $(x,b)\in C_1$ then $(\bar{x},g(x,b))\in C_1$ and $\bar{\bar{x}}=x$}\label{eq:*}\\
\text{if $(y,x,b)\in C_2$ then $\bar{y}(yx)=x,\, (yx)\bar{x}=y,\, x(\overline{yx})=\bar{y},\, (\overline{yx})y=\bar{x}$}\label{eq:**}
\end{eqnarray}
\end{proposition}

Note that when $\bar{x}$ is unique with the properties $\bar{x}(x\bar{x})=\bar{x}$ and $(x\bar{x})x=x$ then $(\bar{x},g(x,b))\in C_1$ as soon as $(x,b)\in C_1$.

\vspace{.5cm}
From now on we assume the conditions of Proposition~\ref{prop:3} and $\bar{1}=1$.

\begin{proposition}\label{prop:4}
If $(1,b)\in C_1$ for every $b\in B$, then the pair $(m\theta,m\varphi)$ is bi-exact. Furthermore, the triple $(\theta,\varphi,m\colon{C_2\to C_1})$ is an associative involutive-2-link as soon as $(X,\cdot)$ is a semigroup. 
\end{proposition}

\begin{proposition}\label{prop:5} If $(1,b)\in C_1$ for all $b\in B$ and the pairs $(m,m\theta)$, $(m,m\varphi)$ are jointly monomorphic, then  the triple $(\theta,\varphi,m\colon C_2\to C_1)$ is a unital involutive-2-link as soon as $(X,\cdot,1)$ is a unital magma.
\end{proposition}

Merging the two previous results, while using Theorem \ref{thm:one}, we obtain.

\begin{proposition}\label{prop:6} If $(1,b)\in C_1$ for all $b\in B$, the pairs $(m,m\theta)$, $(m,m\varphi)$ are jointly monomorphic and $(X,\cdot,1)$ is a monoid, then the triple $(\theta,\varphi,m\colon C_2\to C_1)$ is a unital and associative involutive-2-link in the category of sets and maps.
Moreover, the underlying reflexive graph of the internal groupoid associated with $(\theta,\varphi,m\colon C_2\to C_1)$  is 
\begin{equation}\label{diag: concrete underlying reflexive graph}\xymatrix{C_1   \ar@<1ex>[r]^{d} \ar@<-1ex>[r]_{g} & B \ar[l]|{e},}
\end{equation}
with $d(x,b)=b$ and $e(b)=(1,b)$.
\end{proposition}

 In order to lift the structure $(\theta,\varphi,m\colon C_2\to C_1)$ from the category of sets and maps to the category of magmas and magma homomorphisms we will now assume that $B=(B,+)$ is  a magma and consider the sets $C_1$ and $C_2$ as magmas with operations, respectively,
\begin{eqnarray*}
(x,b)+(x',b')=(f(x,b,x',b'),b+b'),\\
(y,x,b)+(y',x',b')=(f(y,g(x,b),y',g(x',b')),f(x,b,x',b'),b+b').
\end{eqnarray*}

For simplicity, let us from now on assume that $(X,\cdot,1)$ is a group, with  $\Bar{x}=x^{-1}$, and $f(1,b,1,b')=1$ for all $b,b'\in B$. 

\begin{proposition}\label{prop:7} The structure $(\theta,\varphi,m\colon C_2\to C_1)$ is a unital and associative involutive-2-link in the category of magmas if and only if the condition
\begin{equation}\label{eq: g=g+g}
g(f(x,b,x',b'),b+b')=g(x,b)+g(x',b')
\end{equation}
is satisfied for every $(x,b)$ and $(x',b')$ in $C_1$. Furthermore, if $(x,0)\in C_1$ for every $x\in X$, then, for every $(x,b)\in C_1$,   $g(x,b)=g(x,0)+b$.
\end{proposition}


Finally, in order to compare the previous results with the classical notion of crossed module, let us assume that $(B,+,0)$ is a group. It then follows that the map
$f$ is necessarily of the form
\begin{equation}
f(x,b,x',b')=x\cdot f(1,b,x',0)
\end{equation}
and that $\xi(b,x')=f(1,b,x',0)$ is a group action of $B$ on $X$. Moreover, it is not difficult to see that the conditions $(\ref{eq:1 ff}) $ and $(\ref{eq: g=g+g})$ reproduce the classical crossed-module constrains.

\vspace*{1cm}

We have shown that independently of its environment, the category of internal groupoids is equivalent to the full subcategory of  involutive-2-links that are unital and associative. Our approach contrasts with the one in which a groupoid is seen as a reflexive graph equipped with an extra structure, usually adopted when groupoids are studied from an algebraic point of view \cite{BB, Bourn13, Bourn87, Bourn91,Janelidze,NMF.16,NMF.17,NMF.20a}. In our case, the underlying reflexive graph of a groupoid is found as a property of its associated involutive-2-link which  is closer to a more geometrical (or differential) point of view \cite{ Brown1,Brown,Ehresmann,Grothendieck}. However, this work also goes into the  direction of \cite{BournJanel98,Janelidze91} in the sense that it does not require the ambient categories to have pullbacks as canonical constructions and furthermore it can be generalized to higher dimensions \cite{Loday}. 

We conclude with the observation that although an internal groupoid is an instance of an internal category, Brandt \cite{Brandt1926,Brandt} predates Eilenberg and Mac Lane~\cite{Eilenberg-MacLane} in delineating an axiomatic portrait of a (connected) groupoid (\cite{Voight} Remark 19.3.12). Our  approach suggests that internal groupoids can be studied as involutive-2-links in which the unitary and associativity properties, being independent of each other, give rise to a wide spectrum of generalizations.

\section*{Acknowledgement}




  Funded by FCT/MCTES (PIDDAC) through the following Projects: Associate Laboratory ARISE LA/P/0112/2020; UIDP/04044/2020; UI\-DB/04044/2020; PAMI–ROTEIRO/0328/2013 (Nº 022158); MATIS (CENTRO-01-0145-FEDER-000014 - 3362); CENTRO-01-0247-FED\-ER-(069665, 039969); POCI-01-0247-FE\-DER-(069603, 039958, 03\-9863, 024533); Generative thermodynamic; by CDRSP and ESTG from the Polytechnic of Leiria.

  Special thanks are due to IPLeiria's Run-EU program and in particular to the kind and inspiring hospitality offered by FH Vorarlberg -- University of Applied Sciences, at Dornbirn, Austria.

\end{document}